\documentclass[12pt]{article} 

\usepackage{amsmath}
\usepackage{amssymb, amsfonts}

\newtheorem{lemma}{Lemma}[section]
 
\newtheorem{theorem}[lemma]{Theorem}
\newtheorem{corollary}[lemma]{Corollary}

\begin{document}
\title{The Fundamental Theorems in the framework of Bicomplex Topological Modules}
\author{Rajeev Kumar, Romesh Kumar and Dominic Rochon} 
\date{~}
\maketitle

\renewcommand{\thefootnote}{\fnsymbol{footnote}}


\begin{abstract}
In this paper, we generalize the fundamental theorems of functional analysis to the framework of bicomplex topological modules.
\end{abstract}

\section{Introduction}
\setcounter{equation}{0}
This section summarizes a number of known results on the set of bi-complex numbers.

The set $\mathbb{T}$ of bicomplex numbers is defined as 
$$
\mathbb{T} := \{ w= z_1 + \iota_2  z_2: z_1, \; z_2 \in \mathbb{C} (\iota_1)  \},
$$
where $\iota_1$ and $\iota_2$ are independent imaginary units such that $\iota_1^2 = \iota_2^2 = -1.$ The hyperbolic number is denoted as
$
j= \iota_1  \iota_2$ such that $j^2 =1 .$ Under the usual addition and multiplication of bicomplex numbers, $\mathbb{T}$ is a commutative ring. Note that the spaces 
$$
\mathbb{C} (\iota_k) : = \{ a + b  \iota_k ; a, \; b \in \mathbb{R} \}, \, k=1, \, 2
$$ and the set of hyperbolic numbers
$$
\mathbb{D}: = \{ a + d  j; a, \; d \in \mathbb{R} \}.
$$
 are the subrings of the space $\mathbb{T}.$ The norm on $\mathbb{T}$ is defined as 
$$
\mid w \mid := \sqrt{\mid z_1 \mid^2 + \mid z_2 \mid^2} = \sqrt{a^2 + b^2 + c^2 + d^2},
 $$ 
where $w = a + b  \iota_1 + c  \iota_2 + d  \iota_1  \iota_2,$ for $a, \, b, \, c, \, d  \in \mathbb{R}.$
The norm $\mid  . \mid $ is such that 
\[
\mid s \cdot t \mid \le \sqrt{2} \mid s \mid \cdot \mid t \mid
\]

\vspace{3mm} \noindent
{\bf Idempotent basis.} We introduce two bicomplex numbers $e_1$ and $e_2$ defined as 
\[
e_1 = \frac{1 +j }{2} \; {\rm and } \; e_1 = \frac{1 -j }{2}
\] 
Note that these are also hyperbolic numbers such that
\begin{equation} \label{eq:one}
e_k^2 = e_k , \; e_{k}^{+_3}=e_k, \; e_1+ e_2 = 1, \; e_1 .e_2=0, \; {\rm for} \; k= 1, \; 2.
\end{equation}

Any bicomplex number $w$ can be written as 
\[
w= z_1 + z_2 \iota_2= z_{ \hat{1} } \cdot e_1 + z_{ \hat{2} } \cdot e_2,
\]
where
\[
z_{ \hat{1} } = z_1 - z_2 \iota_1 \; {\rm and} \; z_{ \hat{2} }= z_1 + z_2 \iota_1
\]
both are the elements of $\mathbb{C}(\iota_1) .$ We then have another representation of the modulus (norm) on $\mathbb{T}$ as under:
\[
\mid w \mid = \frac{ \sqrt{ \mid z_{ \hat{1} } \mid + \mid z_{ \hat{2} } \mid  }  }{ \sqrt{2} }.
\]
Note that $\mid e_1 \mid= 1/ \sqrt{2} = \mid e_2 \mid.$

\vspace{3mm} \noindent
{\bf Definition 1.} A number $w_1 = z_1 + z_2 \iota_2 \in \mathbb{T}$ is said to have a multiplicative inverse in $\mathbb{T},$ if there exists a number $w_2 = \nu_1 + \nu_2 i \in \mathbb{T}$ such that $w_1  w_2 =1.$ Such elements of $\mathbb{T}$ are called nonsingular elements, otherwise they are called singular. 

By equation (\ref{eq:one}), we see that $\mathbb{T}$ is not even an integral domain so that it is not a field. 

We have two interesting principal ideals $I_1$ and $I_2$ defined as 
\[
I_1= \{ w \cdot e_1 ; \;  w \in \mathbb{T} \}
\]
and 
\[
I_2= \{ w \cdot e_2 ; \; w \in \mathbb{T} \}.
\]
Note that $I_1 \cap I_2 = \{ 0 \}$ and the set $\hat{O}_2$ of all singular elements in $\mathbb{T}$ is nothing but
\[
\hat{O}_2 = I_1 \cup I_2.
\] 
The set set $\hat{O}_2$  is also  caled  Null cone.

Recently,  Lavoie,  Marchildon and  Rochon,  \cite{GMR2}  \cite{GMR1} have introduced bicomplex Hibert spaces and studied some of their basic properties. In this paper  we introduce the concept of bicomplex topologial modules and extend the Principle of uniform boundedness, open mapping theorem, closed graph theorem and Hahn Banach Theorem to this framework. For basic properties of bicomplex analysis  one can refer to  \cite{GMR2} , \cite{GMR1}, \cite{Rochon1},  \cite{T} and references therein.

\section{Principle of Uniform Boundedness}
\setcounter{equation}{0}

\vspace{3mm} \noindent

In this section we study the Uniform Boundedness Principle on bicomplex topological modules.

\vspace{3mm} \noindent
{\bf Definition 2.}
A module $M= (M \, + , \, \cdot)$ over a ring $R$ is called a topological $R$-module if there exists a topology $\tau_M$ on $M$ such that the corresponding operations $+: M \times M \rightarrow M$ and $\cdot : M \rightarrow M $ are continuous.

\vspace{3mm} \noindent {\bf Remark.}
In particularly, if we take $R = \mathbb{T},$ the ring of Bicomplex numbers, we call $M= (M , \, + , \, \cdot)$ as a topological $\mathbb{T}$-module or a topological Bicomplex module. Throughout this paper we assume that $M$ is a topological $\mathbb{T}$ module. With this assumption, we are in the framework of free modules.

\vspace{3mm} \noindent {\bf Example.}
For each $n \in \mathbb{N},$ the set $\mathbb{T}^n := \mathbb{T} \times \mathbb{T} \times \ldots \times \mathbb{T},$ is a topological $\mathbb{T}$-module.

\vspace{3mm} \noindent
{\bf Definition 3.}
Let $M$ be a topological $\mathbb{T}$-module. A map $T: M \rightarrow M$ is said to be a $\mathbb{T}$-linear map or operator if it satisfies the following properties:
\begin{itemize}
\item[$(i)$]  $T(x +y )= T(x)+ T(y), $ for each $x , \; y \in M.$
\item[$(ii)$] $T(\alpha x)= \alpha T(x),$ for each $\alpha \in \mathbb{T}.$
\end{itemize} 

We say that such a $\mathbb{T}$-linear map is bounded if it takes bounded sets into bounded ones.

\vspace{5mm}


The proofs of the following lemmas are straightforward.  
\begin{lemma} \label{th:lone}
For any $a \in M,$ the map $T_a : M \rightarrow M$ defined by $T_a (x) = a + x ,$ for each $x \in M$ is a homeomorphism.
\end{lemma}
\hfill $\Box$

\begin{lemma} \label{th:ltwo}
For each $\lambda \in \mathbb{T},$ the map $M_{\lambda} : M \rightarrow M$ defined by $M_{\lambda} (x) = \lambda \cdot x$ is continuous. In case $\lambda \not \in \hat{O}_2$ then the map $M_{\lambda }$ is a homeomorphism. 
\end{lemma}
\hfill $\Box$

\begin{lemma} \label{th:lthree}
The closure of a submodule of a topological $\mathbb{T}$-module $M$ is a topological $\mathbb{T}$-module.
\end{lemma}
\hfill $\Box$

\noindent
{\bf Remark.} If $M= ( M, \, + , \, \cdot)$ is topological $\mathbb{T}$-module. Then it can be seen that 
\[
V_1 = \{ e_1 w ;  \; w \in M \}
\]
and
\[
V_2 = \{ e_2 w ;  \;  w \in M \}
\]
are two topological vector spaces over the field $F = \mathbb{C} (\iota_1) .$ We see that $M' = V_1 \oplus V_2 $ is also a topological vector space over the field $F = \mathbb{C} (\iota_1).$ In this case note that $M= M' $ as a set.

From Lemma \ref{th:lthree}, we conclude that for any set $B \subseteq M,$ the set $\overline{sp (B)}$ is a topological $\mathbb{T}$-submodule of $M.$ Moreover if 
$\overline{sp (B)}=M,$ we can say that $B$ is a fundamental set.

\vspace{3mm} \noindent
{\bf Definition 4.}
If $V_1$ is spanned by a set $B_1==\{ x_i :  i \in \wedge \}$ in $M$ and $V_2$ is spanned by a set $B_2=\{ y_i :  i \in \wedge \}$ in $M,$  where $\wedge$ is an index set , then the set 
$$
B= B_1 \oplus B_2 :=  \{  x_i + y_i : x_i \in B_1 ,  \; y_i \in B_2 , \; i \in \wedge \}
$$ is said to be a spanning set denoted by 
$sp (B).$   \\
By the Lemma \ref{th:ltwo}, we see that for any set $B \subseteq M, $ the set $\overline{sp (B)}$ is a topological submodule of $M$ over $\mathbb{T}.$ If $\overline{sp (B)}=M,$ we say that $B$ is a fundamental set.

Note that here we have used the fact that $dim (V_1)$ and $dim(V_2)$ is same so that $Card (B_1)= Card (B_2)= Card (B).$ See \cite{GMR1} for details in the finite dimensional case.

\begin{lemma}
The closed topological submodule of $M$ determined by a denumerable set $B$ in a topological $\mathbb{T}$-module is separable.
\end{lemma}
{\bf Proof.}
The set $ A_1 = \{ e_1 b : b \in B \}$ generates a separable topological vector space in $V_1$ and so does $A_2 = \{ e_2 b : b \in B \} $ in $V_2.$ Clearly $B= A_1 \oplus A_2 $ does so in $M.$ 
\hfill $\Box$

\vspace{3mm} \noindent
{\bf Definition 5.}
A set $B$ in a topological $\mathbb{T}$-module $M$ is said to be bounded if given any neighbourhood $V$ of $0$ in $M,$ there exists a number $\epsilon >0$ such that $\alpha B \subseteq M,$ whenever $\mid \alpha \mid \le \epsilon, $ for $\alpha \in \mathbb{T}.$

\begin{lemma} 
A compact subset of a topological $\mathbb{T}$-module $M$ is bounded.
\end{lemma}
{\bf Proof.}
Let $B$ be a compact set in topological $\mathbb{T}$-module $M,$ and let $V$ be a neighbourhood of $0$ in $M.$
Let $B_1 = e_1 B$ and $B_2 = e_2 B.$ Then using \cite[Lemma 8, p-51]{DS}, we see that the sets $B_1$ and $B_2$ are bounded in the topological vector spaces $V_1$ and $V_2,$ respectively. Thus the set $B= B_1 \oplus B_2$ is bounded in $M.$
\hfill $\Box$

\begin{corollary}
A convergent sequence in a topological $\mathbb{T}$-module is bounded.
\end{corollary}
\hfill $\Box$


\vspace{3mm} \noindent
{\bf Definition 6.}
An $F$-module space or a module space of type $F,$ is a topological space $M$ which is also a metric space under some metric $\rho$ such that
\begin{itemize}
\item[$(i)$]  $\rho$ is translation invariant, that is, $\rho (x, \; y) = \rho (x-y, \; 0).$
\item[$(ii)$] $( X, \; \rho)$ is a complete metric space.
\end{itemize} 
In this case, we define an $F$-norm on $X$ as:
\[
\mid x \mid = \rho (x , \; 0).
\]

Note that it is not apparent whether $M$ is a topological $\mathbb{T}$-module until Theorem \ref{th:mm}.

The next theorem demonstrates the principle of uniform boundedness in the setting of topological $\mathbb{T}$-modules.

\begin{theorem}
For each $a \in \wedge,$ where $\wedge$ is an index set,  let $T_a : M \rightarrow M$ be a continuous $\mathbb{T}$-linear map. If for each $x \in M,$ the set $B_x= \{ T_a x : a \in \wedge  \}$ is bounded. Then $\lim_{x \rightarrow 0} T_a x = 0$ uniformly for $a \in \wedge.$
\end{theorem}
{\bf Proof.}
Proof is along the similar lines as in \cite[Theorem 11, p- 52]{DS}.
\hfill $\Box$

\begin{theorem} \label{th:mm}
An $F$-module space is a topological $\mathbb{T}$-module.
\end{theorem}
{\bf Proof.}
Proof is along the similar lines as in \cite[Theorem 12, p- 52]{DS}.
\hfill $\Box$

\begin{theorem} \label{th:am}
A $\mathbb{T}$-linear map of one $F$-module space to another is continuous if and only if it maps bounded sets into bounded ones.
\end{theorem}
{\bf Proof.}
Proof is along the similar lines as in \cite[Theorem 14, p- 52]{DS}.
\hfill $\Box$

\begin{corollary}
Any continuous $\mathbb{T}$-linear map from one topological $\mathbb{T}$-module to another sends  bounded sets into bounded ones.
\end{corollary}
\hfill $\Box$

\begin{corollary}
Any continuous $\mathbb{T}$-linear map from one $F$-module space to another which sends sequences converging to $0$ into bounded sets is continuous.
\end{corollary}
\hfill $\Box$

\begin{theorem} 
Let $(T_{n})_{n=1}^{\infty}$ be a sequence of continuous $\mathbb{T}$-linear maps of one $F$-module space $X$ into another $F$-module space $Y,$ such that the limit
$T(x)= \lim_{n \rightarrow \infty} T_n x $
exists for each $x \in X,$ then $\lim_{x \rightarrow 0} T_n x=0$ uniformly for $n \in \mathbb{N}$ and that $T$ is a continuous $\mathbb{T}$-linear map of $X$ into $Y.$
\end{theorem}
{\bf Proof.}
Proof is along the similar lines as in \cite[Theorem 17, p- 52]{DS}.
\hfill $\Box$

\vspace{3mm} \noindent {\bf Remark.}
The above result also holds in case we replace the sequence $(T_{n})_{n=1}^{\infty}$ by a net $(T_{a})_{a \in \wedge}.$ We state the result for more clarity.

\begin{theorem} 
Let $(T_{a})_{a \in \wedge}$ be a generalized sequence of continuous $\mathbb{T}$-linear maps of one $F$-module space $X$ into another $F$-module space $Y,$ such that the limit
$T(x)= \lim_{n \rightarrow \infty} T_n x $
exists for each $x$ in a fundamental set $B$ in $X,$ and if for each $x \in X$ the set $\{ T_a x\}_{a \in \wedge}$ is bounded, then  the limit 
$T(x)= \lim_{n \rightarrow \infty} T_n x$ exists for each $x \in X,$ and is a continuous $\mathbb{T}$-linear map of $X$ into $Y.$
\end{theorem}
\hfill $\Box$


\section{The Interior Mapping Theorem}
\setcounter{equation}{0}
The interior mapping principle is stated in the following thoeorem.

\begin{theorem}
A continuous $\mathbb{T}$-linear map of one $F$-module space $X$ onto another $F$-module space $Y$ is an open map.
\end{theorem}
{\bf Proof.}
Proof is along similar lines as in \cite[Theorem 1, p-55]{DS}.
\hfill $\Box$

\begin{theorem} \label{th:inv}
A continuous $\mathbb{T}$-linear bijective map of one $F$-module space $X$ onto another $F$-module space $Y$ has a continuous inverse.
\end{theorem}
{\bf Proof.} Proof is obvious by the above arguments or results.
\hfill $\Box$

\vspace{3mm} \noindent
{\bf Note} that the $\mathbb{T}$-linear maps defined by
\[
T_{a} x =ax, \; {\rm for} \; a \in \hat{O}_2
\]
is clearly not onto.


\vspace{3mm} \noindent
{\bf Definition 7.} Let $M$ and $N$ be two topological $\mathbb{T}$-modules.
Let $T$ be a $\mathbb{T}$-linear map  whose domain $D(T)$ defined as
\[
D(T)= \{ x \in M : Tx \in N \}
\]
is a topological $\mathbb{T}$-submodule in $M$ and whose range lies in $N.$ Then the graph of $T$ is the set of all points in $M \times N$ of the form 
$[ x , \; Tx ]$ with $x \in D(T).$ 
The $\mathbb{T}$-linear operator $T$ is said to be closed if its graph is closed in the product space $M \times N.$ An equivalent statement is 
as follows: \\
The $\mathbb{T}$-linear operator $T$ is closed if whenever $x_n \in D(T), \; x_n \longrightarrow x, \; T x_n \longrightarrow y \Rightarrow  x \in D(T) $ and $Tx =y.$ 

Note that the product $M \times N$ of two $F$-module spaces $M$ and $N$ over the same ring $\mathbb{T}$ is also an $F$-module space over $\mathbb{T}$ under the metric $d(\cdot , \, \cdot)$ defined on $M \times N$ as follows: 
\[
d([x , \; y], \; [ x' , \; y' ])= \mid  x - x' \mid + \mid y - y' \mid.
\]


The next result is closed graph theorem in the setting of $F$-module spaces.

\begin{theorem} 
A closed linear map defined on all of an $F$-module space $M$ into an $F$-module space $N$ is continuous.
\end{theorem} 
{\bf Proof.} 
Clearly the graph $G$ of $T$ is a closed $\mathbb{T}$-submodule in the product $F$-module space $M \times N,$ hence $G$ is a complete metric space. Thus $G$ is an $F$-module space. The map $p_M : [x , \; Tx] \mapsto x $ of $G$ onto $M$ is one-to-one, linear, and continuous. Hence, by Theorem \ref{th:inv}, its inverse $p_{M}^{-1}$ is continuous. Thus $T= p_N p_{M}^{-1}$ is continuous by \cite[p-32]{DS}.   
\hfill $\Box$

\begin{theorem}
If a module space $M$ is an $F$-module space under each of the two metrics $\rho_1$ and $\rho_2,$ and if one of the corresponding topologies contains the other, then the two topologies coincide. 
\end{theorem} 
{\bf Proof.} 
Let $\tau_1$ and $\tau_2$ be two metric topologies on the module space $M$ such that $M_1= (M , \; \tau_1)$ and $M_2= (M , \; \tau_2)$ are $F$-module spaces over $\mathbb{T}.$ If $\tau_1 \subseteq \tau_2,$ then the $\mathbb{T}$-linear map $x \mapsto x$ of $M_2$ onto $M_2$ is continuous. By Theorem \ref{th:inv}, it is a homeomorphism so that $\tau_1= \tau_2.$ Hence the theorem.
\hfill $\Box$


\vspace{3mm} \noindent
{\bf Definition 8.}
A family $\mathcal{F}$ of functions which map one module space $M$ into another  module space $N$ over the same ring $\mathbb{T}$ is called total if $f (x)=0, \; \forall f \in \mathcal{F} \Rightarrow  x =0 $ is the only possibility.

\begin{theorem} 
Let $X, \; Y$ and $Z$ be $F$-module spaces over the same ring $\mathbb{T}$ and let $\mathcal{F}$ be a total family of continuous $\mathbb{T}$-linear maps of $X$ into $Y.$ Let $T$ be a linear $\mathbb{T}$-map from $Z$ to $X$ such that $f \circ T$ is continuous $\forall f \in \mathcal{F},$ then $T$ is continuous.
\end{theorem} 
{\bf Proof.} 
Let $w_n \longrightarrow w$ and $T w_n  \longrightarrow x.$ Then 
\[
\lim_{n \rightarrow \infty} f(T w_n )= f(x) \; \forall \; f \in \mathcal{F},
\]
since each $f \circ T$ is continuous so that 
$$
f \circ T (w) = f(x), \; \forall f \in \mathcal{F} \Rightarrow  f(TW-x)=0, \; \forall f \in \mathcal{F} \Rightarrow  Tw=x.
$$
Therefore $\mathcal{F}$ is total.
\hfill $\Box$


 
\section{The Hahn-Banach Theorem}
\setcounter{equation}{0}
First we define $\mathbb{T}$-normed modules as follows:

\vspace{3mm} \noindent
{\bf Definition 9.}
A topological $\mathbb{T}$-module is said to be $\mathbb{T}$-normed module space if there exists a map $\| \cdot \| : M \mapsto \mathbb{R}^+= [0, \; \infty )$ called a $\mathbb{T}$-norm on $M$ if
\begin{itemize}
\item[$(i)$]  $\| \cdot \| : M \rightarrow \mathbb{R}^+ $ is a norm over the field $\mathbb{C}(\iota_1)$ or the field $\mathbb{C}(\iota_2) .$
\item[$(ii)$] $\| w x \| \le \sqrt{2} \mid w \mid \| x \|,$ for each $w \in \mathbb{T}$ and $x \in M.$
\end{itemize} 

Note that $M$ is a topological vector space over the field $\mathbb{C}(\iota_1)$ or the field $\mathbb{C}(\iota_2) .$

A complete $\mathbb{T}$-normed module space is called a Bicomplex Banach module or a $\mathbb{T}$-Banach module.
See \cite{GMR2} and \cite{GMR1} for more light on this aspect.

A $\mathbb{T}$-Banach module is an $F$-module space over $\mathbb{T}$ with the properties as follows:
\[
\| \alpha x \| = \mid \alpha \mid \| x \| , \; \forall \alpha \in \mathbb{C}(\iota_1), \; x \in M
\] 
and
\[ 
\| \alpha x \| \le \sqrt{2} \mid \alpha \mid \| x \| , \; \forall \alpha \in \mathbb{T}, x \in M. 
\]

\begin{lemma}
A set $B$ in a $\mathbb{T}$-normed module space is bounded if and only if $\sup_{ x \in B} \| x \| < \infty.$
\end{lemma}
{\bf Proof.}
A neighbourhood of origin in $M$ contains an $\eta$-neighbourhood 
$$
S_{\eta} = \{ x \in M : \| x \| < \eta \} 
$$
of $0.$ If $a = \sup_{ x \in B} \| x \| < \infty,$ and $\epsilon = \frac{\eta }{ 2 a},$ then $\alpha B \subseteq V$ whenever $\mid \alpha \mid \le \epsilon,$ so that $B$ is bounded.

For $x \in B,$ we have $\| x \| \le a,$ and so
\[
\| \alpha x \| \le \sqrt{2} \mid \alpha \mid \|x \| \le \sqrt{2} \mid \alpha \mid a \le \sqrt{2} \frac{\eta}{2a} a= \frac{\eta}{\sqrt{2}} < \eta. 
\] 

Conversely, if $B$ is bounded, then there exists an $\epsilon >0$ such that
\[
\alpha B \subseteq \delta_1= \{ x \in M : \| x \| < 1  \} \; \forall \; \mid \alpha \mid \le \epsilon. 
\]
For $x \in B, $ we have 
\[
\epsilon \| x \| = \| \epsilon x \| <  1 \Rightarrow \| x \| < \frac{1}{ \epsilon }. 
\]
This proves the result.
\hfill $\Box$


\begin{lemma} \label{th:cont}
For a $\mathbb{T}$-linear map $T$ between $\mathbb{T}$-normed module spaces $X$ and $Y,$  the following properties are equivalent.
\begin{itemize}
\item[$(i)$]  $T $ is a continuous.
\item[$(ii)$] $T $ is a continuous at a point in $M.$
\item[$(iii)$]  $\sup_{\| x \| \le 1} \| Tx \| < \infty.$
\item[$(iv)$] There exists some $M>0,$ such that $\| Tx \| \le \sqrt{2} M \| x \|,$ for each $x \in M.$
\end{itemize} 
\end{lemma}
{\bf Proof.}
(i) $\Rightarrow$ (ii) is obvious.  \\ (ii) $\Rightarrow$ (i) follows from \cite[Lemma 6, p-51]{DS}. Thus (i) $\Leftrightarrow$ (ii). \\
(i)$\Rightarrow$ (iv). If $T$ is continuous at $0 \in M,$ there exists $\epsilon > 0$ such that 
$\| Tx \| < 1 $ if $\| x \| < \epsilon.$ For any $x \neq 0, $ let $ y = \frac{ \epsilon x}{ 2 \| x \|} \in M,$ then 
$$
 \| y \|= ( \frac{\epsilon}{2 \| x \|}) \|x \|= \frac{ \epsilon }{ 2 } < \epsilon,
$$
so that
\[
\frac{\epsilon}{2 \| x \|} \| Tx \|= \| T ( \frac{ \epsilon x}{ 2 \| x \|} ) \| = \| Ty \| < 1
\]
which further implies that
\[
\| Tx \| < \sqrt{2} ( \frac{ \sqrt{2} }{ \epsilon } ) \| x \|. 
\]
This inequality also holds in case $x =0.$ This proves the required implication. 
\\
Clearly (iv) $\Rightarrow$ continuity at $0 \in M.$ Thus we have (iv) $\Rightarrow$ (ii).  
\\
Thus we have (i) $\Leftrightarrow$ (ii) $\Leftrightarrow$ (iv). 
\\
Now we show that (iii) $\Rightarrow$ (iv). \\
If $M = \sup_{ \| x \| \le 1} \| Tx \| < \infty,$ then for any $x \neq 0 ,$ we have
\[
\| Tx \| = \| x \| \| T( \frac{ x }{ \| x \| }  ) \| \le M \| x \| \le \sqrt{2} M \| x \| .
\] 
This inequality also holds in case $x =0.$ This proves the implication.
\hfill $\Box$


\vspace{3mm} \noindent
{\bf Definition 10.}
If $X$ and $Y$ are topological $\mathbb{T}$-modules, we define
\[
B(X, \, Y)= \{ T: X \rightarrow X: T \; {\rm is~a~continuous} \; \mathbb{T}-{\rm linear~map}  \}.
\] 
In case $Y = \mathbb{T},$ we write $X^* = B (X , \; \mathbb{T}).$ Note here we have used the fact that $\mathbb{T}$ is a topological $\mathbb{T}$-module.

Now we define a $\mathbb{T}$-norm on $T \in B(X, \, Y)$ as follows:
\begin{equation} \label{eq:norm}
\| T \| = \| T \|_{ B(X, \, Y)} = \frac{1}{ \sqrt{2} } \sup_{ \| x \| \le 1 } \| Tx \|.
\end{equation}
We have another representation of norm of $T$ when we express $T$ as $T = T_{\hat{1}} e_1 + T_{\hat{2}} e_2$ 
\[
\| T \| = \sqrt{ \frac{ \| T_{\hat{1}} \| + \| T_{\hat{2}} \| }{ 2 } } ,
\]
where each $ T_{\hat{k}}$ maps $X $ into $\mathbb{C} (\iota_1).$ 

If $\| T \| < \infty,$ we say that $T$ is a bounded $\mathbb{T}$-linear operator.

\vspace{3mm}
By lemma \ref{th:cont}, we have the next result as an easy implication.

\begin{theorem}
A $\mathbb{T}$-linear operator $T$ between two $\mathbb{T}$-normed module spaces $X$ and $Y$ is continuous if and only if it is bounded in sense of Definition $10.$
\end{theorem}
\hfill $\Box$

\vspace{3mm}
If $X, \; Y$ and $Z$ are topological $\mathbb{T}$-modules such that $B : X \rightarrow Y$ and $A : Y \rightarrow Z,$ that is, $A$ contains the range of $B,$ then we have
\[
\| AB \| \le \sqrt{2} \| A \| \| B \|.
\]

\begin{theorem}
Let $X$ and $Y$ be two $\mathbb{T}$-Banach modules and $(T_n)_{n \ge 1}$ be a sequence of bounded operators from $X$ into $Y.$ Then the limit
\[
Tx = \lim_{n \rightarrow \infty} T_n x \; {\rm exists~for~each} \; x \in X
\] 
if and only if we have
\begin{itemize}
\item[$(i)$]  $\lim T_n x  $ exists for each $x$ in a fundamental set, and 
\item[$(ii)$] For each $x \in X,$ we have $\sup_{n} \| T_n x\| < \infty.$
\end{itemize} 
\end{theorem}
{\bf Proof.}
Outline: 
\[
\| TX \| = \lim_{n \rightarrow \infty} \| T_n x \| \le \sqrt{2} \lim \inf_{n \rightarrow \infty} \| T_n \| \| x \|
\]
which further implies that
\[
\| T  \| \le \sqrt{2} \lim \inf_{n \rightarrow \infty} \| T_n \|.
\]
\hfill $\Box$

\begin{lemma}
Let $X $ be a $\mathbb{T}$-normed module and $Y$ be a $\mathbb{T}$-Banach module, then $B(X , \; Y )$ is a Banach bicomplex module under the $\mathbb{T}$-norm as defined in the Definition 6.
\end{lemma}
{\bf Proof.}
Clearly $\| T \| =0 $ if and only if $ T=0.$Further,
\[
\| \alpha T \| = \mid \alpha \mid \| T \|, 
\]
when $ \alpha \in \mathbb{C}(\iota_1),$ and 
\[
\| \alpha T \| \le \sqrt{2} \mid \alpha \mid \| T \|, 
\]
for $ \alpha \in \mathbb{T}.$ Since 
\[
\|  (T+ U) (x) \| \le \| Tx \|  + \| Ux \| \le \sqrt{2} ( \| T\| + \| U \| ) \| x \|,
\]
we have
\[
\| T + U \| = \frac{1}{ \sqrt{2} } \sup_{\| x \| \le 1 } \| (T+ U) (x)\| \le \|T \| + \| U \|.
\]
Choose a Cauchy sequence $(T_n)_{n=1}^{\infty}$ in $B(X , \; Y)$ with $\| T_n - T_m \| < \frac{ \epsilon }{ 2 } $ for all $n , \, m \ge n (\epsilon).$ Then $Tx= \lim_{n \rightarrow \infty} T_n (x)$ exists for each $x,$ that is, we have
\[
\| Tx - T_m x \| < \frac{ \epsilon }{ 2 }
\]
and
\[
\| Tx - T_n x \| \le \| Tx -T_m x \| + \sqrt{2} \| T_m  - T_n \| \| x \| < \sqrt{2} ( \frac{ \epsilon }{ 2 \sqrt{2}} +  \frac{ \epsilon }{ 2 } \| x \| ) ,
\]
since the left side of this inequality is independent of $m,$ it is seen by letting $m \rightarrow \infty$ that $\| T - T_n \| \le \sqrt{2} \epsilon $ for $n \ge n (\epsilon)$ so that $\| T \| < \infty,$ and $\| T - T_n \| \longrightarrow 0.$
\hfill $\Box$


\begin{corollary}
The conjugate $\mathbb{T}$-module space $X^*$ of a $\mathbb{T}$-normed module space is a $\mathbb{T}$-Banach module space.
\end{corollary}
\hfill $\Box$

\vspace{1cm}

The next result is clear and and we recall here for the completeness of our presentation. Note that a module over the field $\mathbb{R}$ is obviously a vector space. 

\begin{theorem}\label{th:real}
Let the real function $p$ on module space $X$ over the ring (field) $\mathbb{R}$ satisfy
\[
p (x +y) \le p (x) + p (y), \; p (\alpha x )= \alpha p(x), \; \alpha \ge 0, 
\; x, \; y \in X.
\]
Let $f : Y \rightarrow \mathbb{R}$ be $\mathbb{R}$-linear map on some submodule $Y$ of $X,$ with $f (x) \le p(x), \; \forall x \in Y.$ Then there exists an $\mathbb{R}$-linear map $F : X \rightarrow \mathbb{R}$ such that
\[
F(x) = f(x), \; x \in Y \; {\rm and} \; F(x) \le p (x), \; x \in X.
\]
\end{theorem}
\hfill $\Box$

The next result is the Hahn-Banach Theorem for $\mathbb{T}$-normed modules.
\vspace{5mm}

\begin{theorem}
Let $Y$ be a submodule space of $\mathbb{T}$-normed module space $X.$ Then for each $y \in Y^* ,$ there exists $x^* \in X^*$ with $x^* |_{Y} = y^*$ and
\[
 \| x^* \| = \| y^* \|. 
\] 
\end{theorem}
{\bf Proof.}
If $X$ is a real linear space, then proof is clear by Theorem \ref{th:real} with $p (x)= \| y^* \| \| x \|,$ for $ x \in X,$ and $f= y^*.$ Assume that $X$ is $\mathbb{T}$-normed module. \\
For each $ y \in Y,$ there exists real linear functions $f_1, \; f_2, \; f_3 $ and $f_4$ such that
\[
y^*( y) = (f_1 (y) + \iota_1 f_2 (y)) + \iota_2 (f_3 (y) + \iota_1 f_4 (y)), \; y \in Y.
\] 
Then for $\alpha, \; \beta \in \mathbb{T},$ and $x, \; y \in Y,$ we have
\[
f_1 (\alpha x + \beta y)= \alpha f_1 (x) + \beta f_1 (y),
\]
and
\[
\mid f_1 (y) \mid \le \mid y^* (y) \mid. 
\]
Regarding $X$ as a real module space and so real vector space, we can apply Theorem \ref{th:real}, to get a real linear extension $F_1$ on $X$ by $F_1 |_{Y}=f_1$ and $\| F_1 \| \le \| y^* \|.$

The function $x^*$ on $\mathbb{T}$-module space $X$ is defined by
\[
x^*( x) = F_1 (x) - \iota_1 F_1 ( \iota_1 x) - \iota_2 F_1 ( \iota_2 x) + \iota_1 \iota_2 F_1 ( \iota_1 \iota_2 x).
\] 
Then clearly $x^*$ is additive, and
\begin{eqnarray*}
x^* (\iota_1 x)
& = & F_1 (\iota_1 x) - \iota_1 F_1 (-x) - \iota_2 F_1 (\iota_1 \iota_2 x)+ \iota_1 \iota_2 F_1 (-\iota_2 x) \\
& = & F_1 (\iota_1 x) + \iota_1 F_1 (x) - \iota_2 F_1 (\iota_1 \iota_2 x) - \iota_1 \iota_2 F_1(\iota_2 x)  \\
& = & -\iota_1^2F_1 ( \iota_1 x ) + \iota_1 F_1 (x) + \iota_1^2 \iota_2 F_1 (\iota_1 \iota_2 x) - \iota_1 \iota_2 F_1 (\iota_2 x) \\
& = & \iota_1 [ F_1 (x) - F_1 ( \iota_1 x ) - \iota_2 F(\iota_2 x ) + \iota_1 \iota_2 F_1 (\iota_1 \iota_2 x) ] \\
& = & \iota_1 x^* (x).
\end{eqnarray*}
Similarly we can show that 
\[
x^* (\iota_2 x) = \iota_2 x^* (x)
\]
and that also we have
\begin{eqnarray*}
x^* (\iota_1 \iota_2 x)
& = & F_1 (\iota_1 \iota_2 x) - \iota_1 F_1 (- \iota_2 x) - \iota_2 F_1 (- \iota_1 x)+ \iota_1 \iota_2 F_1 ( x) \\
& = & \iota_1^2 \iota_2^2 F_1 (\iota_1 \iota_2 x) + \iota_1 F_1 (\iota_2 x) + \iota_2 F_1 (\iota_1 x) + \iota_1 \iota_2 F_1( x)  \\
& = & \iota_1 \iota_1 \iota_2 \iota_2  F_1 ( \iota_1 \iota_2 x ) -  \iota_1 \iota_2^2 F_1 (\iota_2 x) - \iota_1^2 \iota_2 F_1 (\iota_1 x) 
+ \iota_1 \iota_2 F_1 ( x) \\
& = & \iota_1 \iota_2 [ F_1 (x) - \iota_1 F_1 ( \iota_1 x ) - \iota_2 F(\iota_2 x ) + \iota_1 \iota_2 F_1 (\iota_1 \iota_2 x) ] \\
& = & \iota_1 \iota_2 x^* (x).
\end{eqnarray*}
Therefore we have
\[
x^* (wx ) = w x^* (x), \; {\rm for~each} \; w \in \mathbb{T}.
\]
This shows that $x^*$ is a $\mathbb{T}$-linear function and that $x^* : X \rightarrow \mathbb{T}.$

Clearly $x^*$ is an extension of $y^* .$ For $y \in Y,$ we have
\begin{eqnarray*}
f_1 (\iota_1 y) + \iota_1 y f_2 (\iota_1 y) +\iota_2 f_3 (\iota_1 y) + \iota_1 \iota_2 f_4 (\iota_1 y)
& = & y^* (\iota_1 y ) = \iota_1 y^* (y) \\
& = & \iota_1 f_1 (y) - f_2 (y)+ \iota_1 \iota_2 f_3 ( \iota_1 y) - \iota_2 f_4 (y)
\end{eqnarray*}
which shows that
\[
f_2 (y) = - f_1 (\iota_1 y),
\]
and similarly we have
\[
f_3 (y) = - f_1 (\iota_2 y),
\]
and
\[
f_4 (y) = f_1 ( \iota_1 \iota_2 y),
\]
and hence that
\[
y^* (y)= f_1 (y) -\iota_1 f_1 (\iota_1 y) -\iota_2 f_1 (\iota_2 y)+ \iota_1 \iota_2 f_1 (\iota_1 \iota_2 y).
\]
Thus $x^*$ is an extension of $y^*.$

Note that we can express $x^*$ as
\[
x^* = ( x^* )_{\hat{1} } e_1 + ( x^* )_{\hat{2}} e_2, 
\]
where each $( x^* )_{\hat{k}}  : X \rightarrow \mathbb{C} (\iota_1)$ for $k=1, \; 2,$ is an extension of the corresponding maps
$( y^* )_{\hat{k}}  : Y \rightarrow \mathbb{C} (\iota_1),$ where  $y^*$ is given by
\[
y^* = ( y^* )_{\hat{1} } e_1 + ( y^* )_{\hat{2}} e_2, 
\]
and
\begin{equation} \label{eq:nn}
\| ( x^* )_{\hat{k}} \| = \| ( y^* )_{\hat{k}} \|
\end{equation}
by the Hahn-Banach theorem for complex vector spaces or the $\mathbb{C} (\iota_1)$-modules. 

Now, using equation \ref{eq:nn}, we have
\begin{eqnarray*}
\| x^* \| 
& = & \sqrt{ \frac{ \| ( x^* )_{\hat{1} } \| + \| ( x^* )_{\hat{2}} \|     }{   2 }} \\
& = & \sqrt{ \frac{ \| ( y^* )_{\hat{1} } \| + \| ( y^* )_{\hat{2}} \|     }{   2 } } \\
& = & \| y^* \|
\end{eqnarray*}
This proves the desired equality and hence the theorem.
\hfill $\Box$

The following results are easy implications of the above interesting result.

\begin{theorem}
Let $Y$ be a $\mathbb{T}$-submodule of the $\mathbb{T}$-normed module $X.$ Let $x \in X$ be such that
\[
\inf_{y \in Y} \| y - x \| =d >0.
\] 
Then there exists a continuous $\mathbb{T}$-linear map $x^* : X \rightarrow \mathbb{T}$ with
\[
x^* (x) =1, \; \mid x^* \mid =\frac{1}{d}, \; x^* (y)=0 , \; \forall y \in Y. 
\] 
\end{theorem}
\hfill $\Box$

\begin{corollary}
Let $x$ be a vector not in the closed $\mathbb{T}$-submodule $Y$ of the $\mathbb{T}$-normed module $X.$ Then there is  
a continuous $\mathbb{T}$-linear map $x^* : X \rightarrow \mathbb{T}$ with
\[
x^* (x) =1, \;  x^* (y)=0 , \; \forall y \in Y. 
\] 
\end{corollary}
\hfill $\Box$

\begin{corollary}
For each $x \neq 0$ in a $\mathbb{T}$-normed module $X,$ there is a continuous $\mathbb{T}$-linear map $x^* : X \rightarrow \mathbb{T}$ with
$\| x^* \| =1$ and $x^* (x)=\| x \| .$ 
\end{corollary}
\hfill $\Box$

\noindent \vspace{5mm}
{\bf Remark.} Note that in this case also, $X^*$ is non-trivial for a non-trivial $\mathbb{T}$-normed module space $X,$ but it is not so in case of $F$-module spaces over the ring $\mathbb{T}.$

\begin{corollary}
For each $x $ in a $\mathbb{T}$-normed module $X,$ 
\[
\| x \| = \sup_{x^* \in S^*} \mid x^* (x) \mid,
\]
where $S^*$ is the closed unit module sphere in the space $X^*$ conjugate to $X.$
\end{corollary}
\hfill $\Box$




%
\vspace{1cm}
\noindent 
Rajeev Kumar, \\ Department of Mathematics,\\
University of Jammu, Jammu, INDIA. \\
\texttt{E-mail: raj1k2@yahoo.co.in} \\
\noindent 

\vspace{1cm}

 \noindent
Romesh Kumar,\\  Department of Mathematics,\\
 University of Jammu, Jammu - 180 006, INDIA. \\
\texttt{E-mails: $romesh{\_}jammu$@yahoo.com}  ,  \texttt{ $romeshmath$@gmail.com}\\

\vspace{5mm}
\noindent 
Dominic Rochon,\\
D\'epartement de math\'ematiques et d'informatique, Universit\'e du Qu\'ebec, \\
Trois-Rivi\`eres, C.P. 500, Qu\'ebec, Canada, G9A 5H7. \\
\texttt{E-mail: dominic.rochon@uqtr.ca} \\

\end{document}